\documentclass[12pt]{article}

\usepackage{amssymb}
\usepackage{amsfonts}

\usepackage{dsfont}

\usepackage{amsmath,amsthm}

\usepackage{graphics}
\usepackage{graphicx}

\DeclareMathOperator{\Per}{Per}

\DeclareMathOperator{\dist}{dist}

\DeclareMathOperator{\Diff}{Diff}

\begin{document}

\sloppy

\righthyphenmin = 2

\newcommand{\al}{\mbox{$\alpha$}}
\newcommand{\si}{\mbox{$\sigma$}}
\newcommand{\de}{\mbox{$\delta$}}
\newcommand{\om}{\mbox{$\omega$}}
\newcommand{\De}{\mbox{$\Delta$}}
\newcommand{\ep}{\varepsilon}
\newcommand{\lam}{\mbox{$\lambda$}}
\newcommand{\La}{\mbox{$\Lambda$}}
\newcommand{\vp}{\mbox{$\varphi$}}
\newcommand{\gam}{\mbox{$\gamma$}}

\newcommand{\ek}{\exp_{p_i}}
\newcommand{\emk}{\exp^{-1}_{p_{i+1}}}
\newcommand{\ekk}{\exp_{p_{i+1}}}

\newcommand{\RR}{\mbox{$\mathds{R}$}}
\newcommand{\Ss}{\mbox{$\textbf{S}$}}
\newcommand{\Z}{\mbox{$\mathds{Z}$}}
\newcommand{\N}{\mbox{$\mathds{N}$}}
\newcommand{\Cc}{\mbox{${\bf C}$}}
\newcommand{\Ff}{\mbox{${\cal F}$}}
\newcommand{\Bb}{\mbox{${\cal B}$}}
\newcommand{\Kk}{\mbox{${\cal K}$}}
\newcommand{\HH}{\mbox{${\cal H}$}}
\newcommand{\LL}{\mbox{${\cal L}$}}

\newcommand{\KK}{\mbox{${\bf K}$}}

\newcommand{\PS}{\mbox{PerSh}}
\newcommand{\LPS}{\mbox{LipPerSh}}

\newcommand{\sref}[1]{(\ref{#1})}

\title{Periodic shadowing and $\Omega$-stability}

\author{A. V.\ Osipov\footnotemark[1],\; S. Yu.\ Pilyugin\footnotemark[1],\;
 and S. B.\ Tikhomirov\footnotemark[2]
\footnotemark[3]}

\date{}

\footnotetext[1] {Faculty of Mathematics and Mechanics, St.\ Petersburg State University, 
                  University av.\ 28,
                  198504, St.\ Petersburg, Russia}

\footnotetext[2] {Department of Mathematics,
National Taiwan University,  No. 1, Section 4, Roosevelt Road, Taipei
106, Taiwan}

\footnotetext[3] {The research of the third author is supported by NSC (Taiwan)
98-2811-M-002-061}

\maketitle

\begin{abstract}

We show that the following three properties of a diffeomorphism $f$
of a smooth closed manifold are equivalent: (i) $f$ belongs to
the $C^1$-interior of the set of diffeomorphisms having
periodic shadowing property; (ii) $f$ has Lipschitz periodic 
shadowing property; (iii) $f$ is $\Omega$-stable. Bibliography: 20 titles.

\end{abstract}
\medskip



\section{Introduction}

The theory of shadowing of approximate trajectories
(pseudotrajectories) of dynamical systems is now a well
developed part of the global theory of dynamical systems
(see, for example, the monographs [1, 2]).

This theory is closely related to the classical theory
of structural stability. It is well known that 
a diffeomorphism has shadowing property in a
neighborhood of a hyberbolic set [3, 4] and a structurally
stable diffeomorpism has shadowing property
on the whole manifold [5 -- 7].
Analyzing the proofs of the first shadowing results by
Anosov [3] and Bowen [4], it is easy to see that, in a
neighborhood of a hyperbolic set, the shadowing property is
Lipschitz (and the same holds in the case of a structurally
stable diffeomorphism, see [1]).

The shadowing property means that, near a
sufficiently precise approximate
trajectory of a dynamical system, there is an exact
trajectory. One can pose a similar question replacing
arbitrary approximate and exact trajectories by periodic
ones (the corresponding property is called periodic
shadowing property, see [8]).

In this paper, we study relations between periodic
shadowing and structural stability (to be more precise,
$\Omega$-stability).

It is easy to give an example of
a diffeomorphism that is not structurally stable but has 
shadowing property (see [9], for example).
Similarly, there exist diffeomorphisms that are not $\Omega$-stable 
but have periodic shadowing property.

Thus, structural stability is not equivalent to shadowing
(and $\Omega$-stability is not equivalent to 
periodic shadowing).

One of possible approaches in the study of relations between
shadowing and structural stability is the passage to
$C^1$-interiors.
At present, it is known that the $C^1$-interior of the set
of diffeomorphisms having shadowing property coincides
with the set of structurally stable diffeomorphisms [10].
Later, a similar result was obtained for orbital shadowing
property (see [11] for details).

In this paper, we show that the $C^1$-interior of the set
of diffeomorphisms having periodic shadowing property coincides
with the set of $\Omega$-stable diffeomorphisms.

We are also interested in the study of the above-mentioned
relations without the passage to $C^1$-interiors.
Let us mention in this context that Abdenur and Diaz
conjectured that a $C^1$-generic diffeomorphism with shadowing
property is structurally stable; they have proved this 
conjecture for so-called tame diffeomorphisms [12].
Recently, it was proved that Lipschitz shadowing and the so-called
variational shadowing are 
equivalent to structural stability [13, 9].

The second main result of this paper states that  Lipschitz periodic shadowing
property is equivalent to $\Omega$-stability.
\medskip

\section{Main results}

Let us pass to exact definitions and statements.

Let $f$ be a diffeomorphism of a smooth closed manifold $M$
with Riemannian metric $\mbox{dist}$. We denote by $Df(x)$ the differential
of $f$ at a point $x\in M$.

Denote by $T_xM$ the tangent space of $M$ at a point $x$;
let $|v|,\;v\in T_xM$, be the norm generated by the metric
$\mbox{dist}$.

As usual, we say that a sequence $\xi=\{x_i\in M,\;i\in\Z\}$
is a $d$-pseudotrajectory of $f$ if
\begin{equation}
\label{0}
\mbox{dist}(f(x_i),x_{i+1})<d,\quad i\in\Z.
\end{equation}

{\bf Definition 1. } We say that $f$ has {\em periodic
shadowing} property if for any positive $\ep$
there exists a positive $d$
such that if $\xi=\{x_i\}$ is a periodic $d$-pseudotrajectory,
then there exists a
periodic point $p$ such that
\begin{equation}
\label{00}
\mbox{dist}(f^i(p),x_i)<\ep,\quad i\in\Z.
\end{equation}

Denote by $\PS$ the set of diffeomorphisms having
periodic shadowing property.

{\bf Definition 2. } We say that $f$ has {\em Lipschitz periodic
shadowing} property if there exist positive constants $\LL,d_0$
such that if $\xi=\{x_i\}$ is a periodic $d$-pseudotrajectory
 with $d\leq d_0$, then there exists a
periodic point $p$ such that
\begin{equation}
\label{00L}
\mbox{dist}(f^i(p),x_i)\leq \LL d,\quad i\in\Z.
\end{equation}

Denote by $\LPS$ the set of diffeomorphisms having
Lipschitz periodic shadowing property.

Denote by $\Omega S$
the set of $\Omega$-stable diffeomorphisms (it is well known that
$f\in\Omega S$ if and only if $f$ satisfies Axiom A and the 
no cycle condition, see, for example, [14]).
Denote by $\Diff ^1(M)$ the space of diffeomorphisms of $M$
with the $C^1$ topology. For a set $P\subset \Diff ^1(M)$
we denote by $\mbox{Int}^1(P)$ its $C^1$-interior.

Let us state our main result.
\medskip

{\bf Theorem. } $\mbox{Int}^1(\PS)=\LPS=\Omega S$.
\medskip

The structure of the paper is as follows. In Sec. 3,
we prove the inclusion $\Omega S\subset\LPS$.
Of course, this inclusion implies that $\Omega S\subset\PS$.
Since the set $\Omega S$ is $C^1$-open, we conclude
that $\Omega S\subset \mbox{Int}^1(\PS)$.
In Sec. 4, 
we prove the inclusion $\mbox{Int}^1(\PS)\subset\Omega S$.
In Sec. 5, 
we prove the inclusion $\LPS\subset\Omega S$.
\medskip

\section{$\Omega S\subset\LPS$}

First we introduce some basic notation. 
Denote by $\Per(f)$ the set of periodic points of $f$ and
by $\Omega(f)$ the nonwandering set of $f$.
Let $N=\sup_{x\in M}\|Df(x)\|$.

Let us 
formulate several auxiliary definitions and statements.

It is well known that if
a diffeomorphism $f$ satisfies Axiom A, then its nonwandering
set can be represented as a disjoint union of 
a finite number of compact sets:
\begin{equation}
\label{spe}
\Omega(f)=\Omega_1\cup\dots\cup\Omega_m,
\end{equation}
where the sets $\Omega_i$ are so-called basic sets (hyperbolic sets
each of which contains a dense positive semi-trajectory).

We say that a diffeomorphism $f$ has Lipschitz 
shadowing property on a set $U$ 
if there exist positive constants $\LL,d_0$
such that if $\xi=\{x_i,\;i\in\Z\}\subset U$ is a $d$-pseudotrajectory
with $d\leq d_0$, then there exists a
point $p\in U$ such that inequalities (\ref{00L}) hold.

We say that a diffeomorphism $f$ is expansive on a set $U$
if there exists a positive number $a$ (expansivity constant)
such that if two trajectories $\{f^i(p):\;i\in\Z\}$ and
$\{f^i(q):\;i\in\Z\}$ belong to $U$ and the inequalities
$$
\mbox{dist}(f^i(p),f^i(q))\leq a,\quad i\in\Z,
$$
hold, then $p=q$.

The following statement is well known (see [1, 14], for example).
\medskip

{\bf Proposition. } {\em If $\Lambda$ is a hyperbolic set,
then there exists a neighborhood $U$ of $\Lambda$ such that
$f$ has Lipschitz shadowing property on $U$ and
is expansive on $U$.}
\medskip

We also need the following two lemmas (see [15]).
\medskip

{\bf Lemma 1. }{\em Let $f$ be a homeomorpism of a compact
metric space $(X,\dist)$. For any neighborhood $U$ of the
nonwandering set $\Omega(f)$ there exist positive numbers
$B,d_1$ such that if $\xi=\{x_i,\;i\in\Z\}$ is a $d$-pseudotrajectory
of $f$ with $d\leq d_1$ and
$$
x_k,x_{k+1},\dots,x_{k+l}\notin U
$$
for some $l>0$ and $k\in\Z$, then $l\leq B$}.
\medskip

Let $\Omega_1,\dots,\Omega_m$ be the basic sets in 
decomposition (\ref{spe}) of the nonwandering set of
an $\Omega$-stable diffeomorphism $f$.
\medskip

{\bf Lemma 2. }{\em Let $U_1,\dots,U_m$ be disjoint
neighborhoods of the basic sets $\Omega_1,\dots,\Omega_m$.
There exist neighborhoods $V_j\subset U_j$ of the
sets $\Omega_j$ and a number $d_2>0$ such that if
$\xi=\{x_i,\;i\in\Z\}$ is a $d$-pseudotrajectory
of $f$ with $d\leq d_2$ such that $x_0\in V_j$ and
$x_t\notin U_j$ for some $j\in\{1,\dots,m\}$ and some $t>0$, then
$x_i\notin V_j$ for $i\geq t$.}
\medskip

{\bf Lemma 3. }{$\Omega S\subset\LPS$.}
\medskip

{\em Proof. } Apply the above proposition and find
disjoint neighborhoods $W_1,\dots,W_m$ of the basic sets $\Omega_1,\dots,\Omega_m$
in decomposition (\ref{spe}) such that (i) 
$f$ has Lipschitz shadowing property on any of $W_j$
with the same constants $\LL,d^*_0$; (ii)
$f$ is expansive on any of $W_j$
with the same expansivity constant $a$.

Find neighborhoods $V_j,U_j$ of $\Omega_j$ (and reduce $d^*_0$, if
necessary) so that the following
properties are fulfilled:

$\bullet$ $V_j\subset U_j\subset W_j,\quad j=1,\dots,m$;

$\bullet$ the statement of Lemma 2 holds for
$V_j$ and $U_j$ with some $d_2>0$;

$\bullet$ the $\LL d^*_0$-neighborhoods of $U_j$
belong to $W_j$.

Apply Lemma 1 to find the corresponding constants
$B,d_1$ for the neighborhood $V_1\cup\dots\cup V_m$
of $\Omega(f)$.

We claim that $f$ has the Lipschitz periodic shadowing property
with constants $\LL,d_0$, where
$$
d_0=\min\left(d^*_0,d_1,d_2,\frac{a}{2\LL}\right).
$$

Take a $\mu$-periodic $d$-pseudotrajectory $\xi=\{x_i,\;i\in\Z\}$
of $f$ with $d\leq d_0$. Lemma 1 implies that there exists
a neighborhood $V_j$ such that $\xi\cap V_j\neq\emptyset$;
shifting indices, we may assume that $x_0\in V_j$.

In this case, $\xi\subset U_j$.
Indeed, if $x_{i_0}\notin U_j$ for some $i_0$, then 
$x_{i_0+k\mu}\notin U_j$ for all $k$.
It follows from Lemma 2 that if $i_0+k\mu>0$,
then $x_i\notin V_j$ for $i\geq i_0+k\mu$,
and we get a contradiction with the periodicity of $\xi$
and the inclusion $x_0\in V_j$.

Thus, there exists a point $p$ such that inequalities
(\ref{00L}) hold. Let us show that $p\in\Per(f)$.
By the choice of $U_j$ and $W_j$, $f^i(p)\in W_j$
for all $i\in\Z$. Let $q=f^\mu(p)$.
Inequalities (\ref{00L}) and the periodicity of $\xi$ imply that
$$
\mbox{dist}(f^i(q),x_{i})=
\mbox{dist}(f^i(q),x_{i+\mu})\leq \LL d,\quad i\in\Z.
$$
Thus,
$$
\mbox{dist}(f^i(q),f^i(p))\leq 2\LL d\leq a,\quad i\in\Z,
$$
which implies that $f^\mu(p)=q=p$. This completes the proof.
\medskip

{\bf Remark. } Thus, we have shown that an $\Omega$-stable
diffeomorphism has periodic shadowing property (and its Lipschitz
variant). It must be noted that it was shown in [16] that there
exist $\Omega$-stable diffeomorphisms that do not have
weak shadowing property (hence, they do not have orbital
and usual shadowing properties, see [11] for details).

\section{ $\mbox{Int}^1(\PS)\subset\Omega S$}

In the proof, we refer to the following well-known statement.
Denote by $\mbox{HP}$ the set of diffeomorphisms $f$
such that every periodic point of $f$ is hyperbolic;
let ${\cal F}=\mbox{Int}^1(\mbox{HP})$.
It is known (see [17, 18]) that the set ${\cal F}$
coincides with the set $\Omega S$ of $\Omega$-stable
diffeomorphisms.

Thus, it suffices for us to prove the following statement.
\medskip

{\bf Lemma 4. } $\mbox{Int}^1(\PS)\subset{\cal F}$.
\medskip

{\em Proof. } In the proof of this lemma, as well as in some
proofs below, we apply the
usual linearization technique based on exponential mapping.

Let $\exp$ be the standard exponential mapping on the 
tangent bundle of $M$ and let $\exp_x$ be the corresponding 
mapping
$$
T_xM\to M.
$$

Let $p$ be a periodic point of $f$; denote $p_i=f^i(p)$
and $A_i=Df(p_i)$.

We introduce the mappings
\begin{equation}
\label{1}
F_i=\emk\circ f\circ\ek: T_{p_i}M\to T_{p_{i+1}}M.
\end{equation}
It follows from the standard properties of the exponential
mapping that $D\exp_x(0)=\mbox{Id}$; hence,
$$
DF_i(0)=A_i.
$$
We can represent
$$
F_i(v)=A_iv+\phi_i(v),
$$
where
$$
\frac{|\phi_i(v)|}{|v|}\to 0\mbox{ as } |v|\to 0.
$$

Denote by $B(r,x)$ the ball in $M$ of radius $r$ centered
at a point $x$ and by $B_T(r,x)$ the ball in $T_xM$ of radius $r$ centered
at the origin.

There exists $r>0$ such that, for any $x\in M$, $\exp_x$ is a diffeomorphism
of $B_T(r,x)$ onto its image, and $\exp^{-1}_x$ is a diffeomorphism
of $B(r,x)$ onto its image. In addition, we may assume that $r$ has
the following property.

If $v,w\in B_T(r,x)$, then
$$
\frac{\mbox{dist}(\exp_x(v),\exp_x(w))}{|v-w|}\leq 2;
$$
if $y,z\in B(r,x)$, then
$$
\frac{|\exp^{-1}_x(y)-\exp^{-1}_x(z)|}{\mbox{dist}(y,z)}\leq 2.
$$

Every time, constructing periodic $d$-pseudotrajectories of $f$,
we take $d$ so small that the considered points of our pseudotrajectories,
points of shadowing trajectories, their ``lifts" to tangent spaces,
etc belong to the corresponding balls $B(r,p_i)$ and $B_T(r,p_i)$
(and we do not repeat this condition on the smallness of $d$).

To prove Lemma 4, it is enough for us to show that
$\mbox{Int}^1(\PS)\subset\mbox{HP}$ and to note that the
left-hand side of this inclusion is $C^1$-open.

To get a contradiction, let us assume that a diffeomorphism
$f\in\mbox{Int}^1(\PS)$ has a nonhyperbolic periodic point $p$.
Fix a $C^1$-neighborhood ${\cal N}\subset\PS$ of $f$.

For simplicity, let us assume that $p$ is a fixed point and
that the matrix $A_0=Df(p)$ has an eigenvalue $\lam=1$
(the remaining cases are considered using a similar
reasoning, see, for example, [19]).

In our case, an analog of mapping (\ref{1}), 
$$
F=\exp_p^{-1}\circ f\circ\exp_p: T_{p}M\to T_{p}M,
$$
has the form
$$
F(v)=A_0v+\phi(v).
$$
Clearly, we can find a number $a\in(0,r)$ (recall that
the number $r$ was fixed above when properties of the
exponential mapping were described),
coordinates $v=(u,w)$ in $T_pM$ with one-dimensional $u$, and a diffeomorphism
$h\in{\cal N}$ such that if
$$
H=\exp_p^{-1}\circ h\circ\exp_p
$$
and $|v|\leq a$, then
$$
H(v)=Av=(u,Bw),
$$
where $B$ is a matrix of size $(n-1)\times(n-1)$
(and $n$ is the dimension of $M$).
For this purpose, we take a matrix $A$, close to $A_0$
and having an eigenvalue $\lam=1$ of multiplicity one,
and ``annihilate" the $C^1$-small term
$(A_0-A)v+\phi(v)$ in the small ball $B_T(a,p)$.

Take a positive $\ep$ such that $8\ep<a$.
Since $h\in{\cal N}$, there exists a corresponding $d\in(0,\ep)$
from the definition of periodic shadowing (for the
diffeomorphism $h$). Take a natural number $K$
such that $Kd>8\ep$. Reducing $d$, if necessary, we may assume that
\begin{equation}
\label{2.01}
8\ep<Kd<2a.
\end{equation}
Let us construct a sequence $y_k\in T_pM,\;k\in\Z,$ as follows:
$$
y_0=0,\quad y_{k+1}=Ay_k+\left(\frac{d}{2},0\right),\quad 0\leq k\leq K-1,
$$
$$
y_{k+1}=Ay_k-\left(\frac{d}{2},0\right),\quad K\leq k\leq 2K-1,
$$
and $y_{k+2K}=y_k,\;k\in\Z$. Clearly,
\begin{equation}
\label{2.2}
y_K=\left(\frac{Kd}{2},0\right).
\end{equation}
Let
$$
x_k=\exp_p(y_k).
$$
Since
$$
\exp_p^{-1}(h(x_k))=H(y_k)=Ay_k
$$
and
$$
|y_{k+1}-Ay_k|=\frac{d}{2},
$$
the sequence $\xi=\{x_k\}$ is a $2K$-periodic $d$-pseudotrajectory
of $h$.

By our assumption, there exists a periodic point $p_0$ of $h$
such that
$$
\mbox{dist}(p_k,x_k)<\ep,\quad k\in\Z,
$$
where $p_k=h^k(p_0)$. Let
$$
p_k=\exp_p(q_k),\quad k\in\Z,
$$
where $q_k=(U_k,W_k)$, and let $y_k=(u_k,w_k)$; then 
$$
|U_k-u_k|\leq|q_k-y_k|<2\ep,\quad k\in\Z,
$$
which implies that
$$
|U_0|\leq|q_0|<2\ep.
$$
Since $q_{k+1}=H(q_k)$, $U_k=U_0$ for all $k$ due to the structure of $H$.
We conclude that $|U_K|<2\ep$ and get a contradiction with 
the inequalities $|U_K-u_K|<2\ep$, (\ref{2.01}), and (\ref{2.2}).
The lemma is proved.

\section{ $\LPS\subset\Omega S$}

In this section, we assume that $f\in\LPS$
(with constants $\LL\geq 1,d_0>0$). Clearly, in this case
$f^{-1}\in\LPS$ as well (and we assume that the
constants $\LL,d_0$ are the same for $f$ and $f^{-1}$).

In the construction of pseudotrajectories, we apply the
same linearization technique as in the previous section.
\medskip

{\bf Lemma 5. } {\em Every point $p\in\Per(f)$ is hyperbolic.}
\medskip

{\em Proof. } To get a contradiction, let us assume that $f$ has
a nonhyperbolic periodic point $p$ (to simplify notation, we
assume that $p$ is a fixed point; literally the same reasoning
can be applied to a periodic point of period $m>1$).

In this case, mapping (\ref{1}) takes the form
$$
F(v)=\exp^{-1}_p\circ f\circ\exp_p(v)=Av+\phi(v),
$$
where $A$ is a nonhyperbolic matrix.
The following two cases are possible:

(Case 1): $A$ has a real eigenvalue $\lam$ with $|\lam|=1$;

(Case 2): $A$ has a complex eigenvalue $\lam$ with $|\lam|=1$.

We treat in detail only Case 1; we give a comment concerning
Case 2. To simplify presentation, we assume that 1 is an
eigenvalue of $A$; the case of eigenvalue $-1$ 
is treated similarly.

We can find coordinates $v$ in $T_pM$ such that, with respect to
this coordinate, the matrix $A$
has block-diagonal form,
\begin{equation}
\label{bform}
A=\mbox{diag}(B,P),
\end{equation}
where $B$ is a Jordan block of size $l\times l$: 
$$
B=\left(
\begin{array}{ccccc}
1&1&0&\ldots&0\\
0&1&1&\ldots&0\\
\vdots&\vdots&\vdots&\ddots&\vdots\\
0&0&0&\ldots&1
\end{array}
\right).
$$

Of course, introducing new coordinates, we have to change
the constants $\LL,d_0,N$; we denote the new constants by the
same symbols. In addition, we assume that $\LL$ is integer.

We start considering the case $l=2$; in this case,
$$
B=\left(
\begin{array}{cc}
1&1\\
0&1
\end{array}
\right).
$$
Let
$$
e_1=(1,0,0,\dots,0)
\mbox{ and }
e_2=(0,1,0,\dots,0)
$$
be the first two vectors of the standard orthonormal basis.

Let $K=25\LL$.

Take a small $d>0$ and construct a finite sequence $y_0,\dots,y_Q$
in $T_pM$ (where $Q$ is determined later) as follows:
$y_0=0$ and 
\begin{equation}
\label{pst}
y_{k+1}=Ay_k+de_2,\quad k=0,\dots, K-1.
\end{equation}
Then 
$$
y_K=(Z_1(K)d,Kd,0,\dots,0),
$$
where the natural number $Z_1(K)$ is determined by $K$ (we do not 
write $Z_1(K)$ explicitly). Now we set
$$
y_{k+1}=Ay_k-de_2,\quad k=K,\dots, 2K-1.
$$
Then 
$$
y_{2K}=(Z_2(K)d,0,0,\dots,0),
$$
where the natural number $Z_2(K)$ is determined by $K$ as well.
Take $Q=2K+Z_2(K)$; if we set
$$
y_{k+1}=Ay_k-de_1,\quad k=2K,\dots, Q-1,
$$
then $y_Q=0$.
Let us note that both numbers $Q$ and
$$
Y:=\frac{\max_{0\leq k\leq Q-1}|y_k|}{d}
$$
are determined by $K$ (and hence, by $\LL$).

Now we construct a $Q$-periodic sequence $y_k,k\in\Z,$
that coincides with the above sequence for $k=0,\dots,Q$.

We set $x_k=\exp_p(y_k)$ and claim that if $d$ is small
enough, then $\xi=\{x_k\}$ is a $4d$-pseudotrajectory of $f$
(and this pseudotrajectory is $Q$-periodic by construction).

Indeed, we know that
$|y_k|\leq Yd$ for $k\in\Z$. Since $\phi(v)=o(|v|)$ as $|v|\to 0$,
\begin{equation}
\label{5}
|\phi(y_k)|<d,\quad k\in\Z,
\end{equation}
if $d$ is small enough.

The definition of $\{y_k\}$ 
implies that
\begin{equation}
\label{6}
|y_{k+1}-Ay_{k}|=d,\quad k\in\Z.
\end{equation}

Note that
$$
\exp^{-1}_p(f(x_k))=F(y_k)=Ay_k+\phi(y_k);
$$
thus, it follows from (\ref{5}) and (\ref{6}) that
$$
|y_{k+1}-\exp^{-1}_p(f(x_k))|\leq |y_{k+1}-Ay_{k}|+|\phi(y_k)|<2d,
$$
which implies that  
$\xi=\{x_k\}$ is a $4d$-pseudotrajectory of $f$
if $d$ is small
enough.

Now we estimate the distances between points of trajectories
of the mapping $F$ and its linearization.

Let us take a vector $q_0\in T_pM$ and assume that the
sequence $q_k=F^k(q_0)$ belongs to
the ball $|v|\leq (Y+8\LL)d$ for $0\leq k\leq K$.
Let $r_k=A^kq_0$ (we impose no conditions on $r_k$ since
below we estimate $\phi$ at points $q_k$ only).

Take a small number $\mu\in(0,1)$ (to be chosen later) and assume
that $d$ is small enough, so that the inequality
$$
|\phi(v)|\leq\mu|v|
$$
holds for $|v|\leq (Y+8\LL)d$. 

Then
$$
|q_1|\leq|Aq_0|+|\phi(q_0)|\leq (N+1)|q_0|,\dots,
|q_{k}|\leq|Aq_{k-1}|+|\phi(q_{k-1})|\leq (N+1)^k|q_0|
$$
for $1\leq k\leq K$, and
$$
|q_1-r_1|=|Aq_0+\phi(q_0)-Aq_0|\leq\mu|q_0|,
$$
$$
|q_2-r_2|=|Aq_1+\phi(q_1)-Ar_1|\leq N|q_1-r_1|+\mu|q_1|
\leq \mu(2N+1)|q_0|,
$$
$$
|q_3-r_3|\leq N|q_2-r_2|+\mu|q_2|
\leq \mu(N(2N+1)+(N+1)^2)|q_0|,
$$
and so on.

Thus, there exists a number $\nu=\nu(K,N)$ such that
$$
|q_k-r_k|\leq \mu\nu|q_0|,\quad 0\leq k\leq K.
$$
We take $\mu=1/\nu$, note that $\mu=\mu(K,N)$, and get
the inequalities
\begin{equation}
\label{7}
|q_k-r_k|\leq |q_0|,\quad 0\leq k\leq K,
\end{equation}
for $d$ small enough.

Since $f\in\LPS$, for $d$ small enough, the $Q$-periodic 
$4d$-pseudotrajectory $\xi$ is $4\LL d$-shadowed
by a periodic trajectory. Let $p_0$ be a point of
this trajectory such that  
\begin{equation}
\label{8}
\mbox{dist}(p_k,x_k)\leq 4\LL d,\quad k\in\Z,
\end{equation}
where $p_k=f^k(p_0)$.
Let $q_k=\exp^{-1}_p(p_k)$.

The inequalities $|y_k|\leq Yd$ and (\ref{8}) imply that
\begin{equation}
\label{9}
|q_k|\leq |y_k|+2\mbox{dist}(p_k,x_k)\leq (Y+8\LL)d,\quad k\in\Z.
\end{equation}
Note that $|q_0|\leq 8\LL d$.

Set $r_k=A^kq_0$;
we deduce from estimate (\ref{7}) that if $d$ is small enough,
then
\begin{equation}
\label{10}
|q_K-r_K|\leq |q_0|\leq 8\LL d.
\end{equation}

Denote by $v^{(2)}$ the second coordinate of a vector $v\in T_pM$.

It follows from the structure of the matrix $A$ that
\begin{equation}
\label{11}
|r_K^{(2)}|=|q_0^{(2)}|\leq 8\LL d.
\end{equation}
The relations
$$
|y_K^{(2)}|=Kd\mbox{ and } |q_K-y_K|\leq 8\LL d
$$
imply that
\begin{equation}
\label{12}
|q_K^{(2)}|\geq Kd-8\LL d=17\LL d
\end{equation}
(recall that $K=25\LL$).

Estimates (\ref{10})--(\ref{12}) are contradictory.
Our lemma is proved in Case 1 for $l=2$.

If $l=1$, then the proof is simpler; the first coordinate
of $A^kv$ equals the first coordinate of $v$, and we
construct the periodic pseudotrajectory perturbing
the first coordinate only.

If $l>2$, the reasoning is parallel to that above;
we first perturb the $l$th coordinate to make it $Kd$,
and then produce a periodic sequence consequently making zero 
the $l$th coordinate, the $(l-1)$st coordinate, and
so on.

If $\lam$ is a complex eigenvalue, $\lam=a+bi$, 
we take
a real $2\times 2$ matrix
$$
R=\left(
\begin{array}{cc}
a&-b\\
b&a\\
\end{array}
\right)
$$
and assume that in representation (\ref{bform}),
$B$ is a real $2l\times 2l$ Jordan block:
$$
B=\left(
\begin{array}{ccccc}
R&E_2&0&\ldots&0\\
0&R&E_2&\ldots&0\\
\vdots&\vdots&\vdots&\ddots&\vdots\\
0&0&0&\ldots&R
\end{array}
\right),
$$
where $E_2$ is the $2\times 2$ unit matrix.

After that, almost the same reasoning works;
we note that $|Rv|=|v|$ for any 2-dimensional vector $v$
and construct periodic pseudotrajectories
replacing, for example, formulas (\ref{pst}) 
by the formulas
$$
y_{k+1}=Ay_k+dw_k,\quad k=0,\dots,K-1,
$$
where $j$th coordinates of the vector $w_k$
are zero for $j=1,\dots,2l-2,2l+1,\dots,n$,
while the 2-dimensional vector corresponding
to $(2l-1)$st and $2l$th coordinates has
the form $R^kw$ with $|w|=1$, and so on.
We leave details to the reader. The lemma is proved.
\medskip

{\bf Lemma 6. }{\em There exist constants $C>0$ and $\lam\in(0,1)$ 
depending only on $N$ and $\LL$ and such that,
for any point $p\in\Per(f)$, there exist complementary
subspaces $S(p)$ and $U(p)$ of the tangent space $T_pM$
that are $Df$-invariant, i.e.,

(H1) $Df(p)S(p)=S(f(p))$ and $Df(p)U(p)=U(f(p))$,

\noindent and the inequalities

(H2.1) $|Df^j(p)v|\leq C\lam^j|v|, \quad v\in S(p), j\geq 0$,

\noindent and 

(H2.2) $|Df^{-j}(p)v|\leq C\lam^j|v|, \quad v\in U(p), j\geq 0$,

\noindent hold}.
\medskip

{\bf Remark. } Lemma 6 means that the set $\Per(f)$ has all 
the standard properties of a hyperbolic set, with
the exception of compactness.
\medskip

{\em Proof. } Take a periodic point $p\in\Per(f)$; let $m$
be the minimal period of $p$.

Denote $p_i=f^i(p)$, $A_i = D f(p_i)$, and $B = D f^m(p)$. 
It follows from Lemma 5
that the matrix $B$ is hyperbolic. Denote by $S(p)$ and $U(p)$
the invariant subspaces
of $B$ corresponding to parts of its spectrum inside and outside
the unit disk, respectively. Clearly, $S(p)$ and $U(p)$ are
invariant with respect to $Df$,
$T_{p}M = S(p) \oplus U(p)$, and the following relations hold:
\begin{equation}\label{1.1}
\lim_{n \to +\infty} B^n v_s = \lim_{n \to +\infty} B^{-n} v_u = 0,
\quad v_s \in S(p), v_u \in U(p).
\end{equation}

We prove that inequalities (H2.2) hold with $C=16\LL$ and
$\lam=1+1/(8\LL)$ (inequalities (H2.1) are established by 
similar reasoning applied to $f^{-1}$ instead of $f$).

Consider an arbitrary nonzero vector $v_u \in U(p)$ and an integer $j\geq 0$.
Define sequences $v_i,
e_i \in T_{p_i}M$ and $\lam_i > 0$ for $i\geq 0$ as follows:
$$
v_0 = v_u, \quad v_{i+1} = A_i v_i, \quad e_i = \frac{v_i}{|v_i|},
\quad \lam_i = \frac{|v_{i+1}|}{|v_i|} = |A_i e_i|.
$$
Let
$$
\tau= \frac{\lam_{m-1}\cdot \ldots \cdot \lam_1 + \lam_{m-1}\cdot \ldots \cdot
\lam_2 + \ldots + \lam_{m-1} + 1}{\lam_{m-1}\cdot \ldots \cdot \lam_0}.
$$
Consider the sequence $\{a_i \in \RR,\;i\geq 0\}$ defined by the following
formulas:
\begin{equation}
\label{1.2}
a_0 = \tau, \quad a_{i+1} = \lam_i a_i -1.
\end{equation}
Note that
\begin{equation}
\label{1.3}
a_{m} = 0 \quad \mbox{and} \quad a_i >0, \quad i \in [0, m-1].
\end{equation}
Indeed, if $a_i\leq 0$  for some $i \in [0, m-1]$, then $a_k<0$
for $k \in [i+1, m]$.

It follows from (\ref{1.1}) that there exists $n > 0$ such that
\begin{equation}
\label{2.1}
|B^{-n}\tau e_0| < 1.
\end{equation}

Consider the finite sequence $\{w_i \in T_{p_i}M,\;i\in[0,m(n+1)]\}$
defined as follows:
$$
\begin{cases} 
w_i=a_i e_i, & \quad i \in [0, m-1], \\
w_{m} = B^{-n}\tau e_0, & \;\\
w_{m+1+i} = A_i w_{m+i}, & \quad i \in [0, mn - 1].
\end{cases}
$$
Clearly,
$$
w_{km}=B^{k-1-n}\tau e_0,\quad k\in[1,n+1],
$$
which means that we can consider $\{w_i\}$ as an
$m(n+1)$-periodic sequence defined for $i\in\Z$.

Let us note that
$$
A_iw_i=a_iA_ie_i=
a_i\frac{v_{i+1}}{|v_{i}|},\quad i\in[0,m-2],
$$
$$
w_{i+1}=(\lam_ia_i-1)\frac{v_{i+1}}{|v_{i+1}|}
=a_i\frac{v_{i+1}}{|v_{i}|}-e_{i+1},\quad i\in[0,m-2],
$$
and 
$$
A_{m-1}w_{m-1}=a_{m-1}\frac{v_{m}}{|v_{m-1}|}=
\frac{v_{m}}{\lam_{m-1}|v_{m-1}|}=e_m
$$
(in the last relation we take into account that
$a_{m-1}\lam_{m-1}=1$ since $a_m=0$).

The above relations and condition (\ref{2.1}) imply that
\begin{equation}
\label{15}
|w_{i+1} - A_i w_i| < 2, \quad i \in \Z.
\end{equation}

Now we take a small $d>0$ and consider 
the $m(n+1)$-periodic sequence $\xi=\{x_i=\mbox{exp}_{p_i}(dw_i),\;i\in \Z\}$.

We claim that if $d$ is small enough, then $\xi$
is a $4d$-pseudotrajectory of $f$.

Denote
$$
\zeta_{i+1}=\emk(f(x_i))\;\mbox{ and }\;\zeta'_{i+1}=\emk(x_{i+1}).
$$
Then
$$
\zeta_{i+1}=\emk f(\ek(dw_i))=F_i(dw_i)=A_idw_i+\phi_i(dw_i),
$$
where the mapping $F_i$ is defined in (\ref{1}) and $\phi_i(v)=o(|v|)$, and
$$
\zeta'_{i+1}=\emk(x_{i+1})=dw_{i+1}.
$$
It follows from estimates (\ref{15}) that
$$
|\zeta'_{i+1}-\zeta_{i+1}|\leq 2d
$$
for small $d$, and
$$
\mbox{dist}(f(x_i),x_{i+1})\leq 4d.
$$

By Lemma 5, the $m$-periodic trajectory $\{p_i\}$ is
hyperbolic; hence, $\{p_i\}$ has a neighborhood in which
$\{p_i\}$ is a unique periodic trajectory.
It follows that if $d$ is small enough, then
the pseudotrajectory $\{x_i\}$ is
$4\LL d$-shadowed by $\{p_i\}$.

The inequalities $\dist(x_i,p_i)\leq 4\LL d$
imply that $|a_i|=|w_i|\leq 8\LL$ for $0\leq i\leq m-1$.

Now the equalities $\lam_i =
(a_{i+1}+1)/a_i$ imply that
if $0\leq i\leq m-1$, then
$$
\lam_0\cdot\ldots\cdot\lam_{i-1}
=\frac{a_{1}+1}{a_0}\frac{a_{2}+1}{a_{1}}\dots
\frac{a_{i}+1}{a_{i-1}}=
$$
$$
=\frac{a_{i}+1}{a_0}\left(1+\frac{1}{a_{1}}\right)\dots\left(1+\frac{1}{a_{i-1}}\right)\geq
$$
$$
\geq \frac{1}{8\LL}\left(1+\frac{1}{8\LL}\right)^{i-1}>
\frac{1}{16\LL}\left(1+\frac{1}{8\LL}\right)^{i}
$$
(we take into account that $1+1/(8\LL)<2$
since $\LL\geq 1$).

It remains to note that 
$$
|Df^i(p)v_u|=\lam_{i-1}\cdots\lam_0|v_u|,\quad 0\leq i\leq m-1,
$$
and that we started with an arbitrary vector $v_u\in U(p)$.

This proves our statement for $j\leq m-1$. If $j\geq m$, we take an integer $k>0$
such that $km>j$ and repeat the above reasoning for the periodic trajectory
$p_0,\dots,p_{km-1}$ (note that we have not used the condition that $m$
is the minimal period).
Lemma 6 is proved.
\medskip

\medskip

{\bf Lemma 7. } {\em If} $f\in\LPS$, {\em then $f$ satisfies Axiom A.}
\medskip

{\em Proof. } Denote by $P_l$ the set of points $p\in\Per(f)$
of index $l$ (as usual, the index of a hyperbolic periodic point
is the dimension of its unstable manifold).

Let $R_l$ be the closure of $P_l$. Clearly, $R_l$ is a compact
$f$-invariant set. We claim that any $R_l$ is a hyperbolic set.
Let $n=\mbox{dim}M$.

Consider a point $q\in R_l$ and fix a sequence of points $p_m\in P_l$
such that $p_m\to q$ as $m\to\infty$. By Lemma 6, there exist 
complementary subspaces $S(p_m)$ and $U(p_m)$ of $T_{p_{m}}M$
(of dimensions $n-l$ and $l$, respectively)
for which estimates (H2.1) and (H2.2) hold.

Standard reasoning shows that, introducing local coordinates
in a neighborhood of $(q,T_qM)$ in the tangent bundle of $M$,
we can select a subsequence $p_{m_k}$ for which the 
sequences $S(p_{m_k})$ and $U(p_{m_k})$ converge (in the
Grassmann topology) to subspaces of $T_qM$
(let $S_0$ and $U_0$ be the corresponding limit subspaces).

The limit subspaces $S_0$ and $U_0$ are complementary in $T_qM$.
Indeed, consider the ``angle" $\beta_{m_k}$ between the subspaces
$S(p_{m_k})$ and $U(p_{m_k})$ which is defined (with respect
to the introduced local coordinates in a neighborhood of $(q,T_qM)$)
as follows:
$$
\beta_{m_k}=\min |v^s-v^u|,
$$ 
where the minimum is taken over all possible pairs of
unit vectors $v^s\in S(p_{m_k})$ and $v^u\in U(p_{m_k})$.

It is shown in [16, Lemma 12.1] that the values $\beta_{m_k}$
are estimated from below by a positive constant
$\alpha=\alpha(C,\lam,N)$. Clearly, this implies that the
subspaces $S_0$ and $U_0$ are complementary.

It is easy to show that the limit subspaces $S_0$ and $U_0$
are unique (which means, of course, that the sequences
$S(p_m)$ and $U(p_m)$ converge). For the convenience of the reader,
we prove this statement (our reasoning is close to that of [16]).

To get a contradiction, assume that there is a subsequence
$p_{m_i}$ for which the sequences $S(p_{m_i})$ and $U(p_{m_i})$ converge
to complementary subspaces $S_1$ and $U_1$ different from
$S_0$ and $U_0$ (for definiteness, we assume that $S_0\setminus S_1\neq\emptyset$).

Due to the continuity of $Df$, the inequalities
$$
|Df^j(q)v|\leq C\lam^j|v|,\quad v\in S_0\cup S_1,
$$
and
$$
|Df^j(q)v|\geq C^{-1}\lam^{-j}|v|,\quad v\in U_0\cup U_1,
$$
hold for $j\geq 0$.

Since
$$
T_qM=S_0\oplus U_0=S_1\oplus U_1,
$$
our assumption implies that there is a vector $v\in S_0$
such that 
$$
v=v^s+v^u,\quad v^s\in S_1, v^u\in U_1, v^u\neq 0.
$$
Then
$$
|Df^j(q)v|\leq C\lam^j|v|\to 0,\quad j\to\infty,
$$
and
$$
|Df^j(q)v|\geq C^{-1}\lam^{-j}|v^u|-C\lam^j|v^s|\to \infty,\quad j\to\infty,
$$
and we get the desired contradiction.

It follows that there are uniquely defined complementary subspaces 
$S(q)$ and $U(q)$ for $q\in R_l$ with proper hyperbolity
estimates; the $Df$-invariance of these subspaces is obvious.
We have shown that each $R_l$ is a hyperbolic set with
$\mbox{dim}S(q)=n-l$ and $\mbox{dim}U(q)=l$ for $q\in R_l$.

If $r\in\Omega(f)$, then there exists a sequence of points $r_m\to r$
as $m\to\infty$ and a sequence of indices $k_m\to\infty$
as $m\to\infty$ such that $f^{k_m}(r_m)\to r$.

Clearly, if we continue the sequence 
$$
r_m,f(r_m),\dots,f^{k_m-1}(r_m)
$$
periodically with period $k_m$, we get a periodic
$d_m$-pseudotrajectory of $f$ with $d_m\to 0$ as $m\to\infty$.

Since $f\in \LPS$, for large $m$ there exist periodic points
$p_m$ such that $\mbox{dist}(p_m,r_m)\to 0$ as $m\to\infty$.
Thus, periodic points are dense in $\Omega(f)$.

Since hyperbolic sets with different dimensions of the subspaces
$U(q)$ are disjoint, we get the equality
$$
\Omega(f)=R_0\cup\dots\cup R_{n},
$$
which implies that $\Omega(f)$ is hyperbolic.
The lemma is proved.
\medskip

It was mentioned above that if
a diffeomorphism $f$ satisfies Axiom A, then its nonwandering
set can be represented as a disjoint union of 
a finite number of basic sets (see representation (\ref{spe})).

The basic sets $\Omega_i$ have
stable and unstable ``manifolds":
$$
W^s(\Omega_i)=\{x\in M:\;\mbox{dist}(f^k(x),\Omega_i)\to 0,\quad k\to\infty\}
$$
and
$$
W^u(\Omega_i)=\{x\in M:\;\mbox{dist}(f^k(x),\Omega_i)\to 0,\quad k\to-\infty\}.
$$
If $\Omega_i$ and $\Omega_j$ are basic sets, we write $\Omega_i\to\Omega_j$
if the intersection
$$
W^u(\Omega_i)\cap W^s(\Omega_j)
$$
contains a wandering point.

We say that $f$ has a 1-cycle if there is a basic set $\Omega_i$ such that
$\Omega_i\to\Omega_i$.

We say that $f$ has a $t$-cycle if there are $t>1$ basic sets 
$$
\Omega_{i_1},\dots,\Omega_{i_t}
$$
such that
$$
\Omega_{i_1}\to\dots\to\Omega_{i_t}\to\Omega_{i_1}.
$$
\medskip

{\bf Lemma 8. } {\em If } $f\in\LPS$, {\em then $f$ has no cycles.}
\medskip

{\em Proof. } To simplify presentation, we prove that $f$ has
no 1-cycles (in the general case, the idea is literally the
same, but the notation is heavy).

To get a contradiction, assume that
$$
p\in(W^u(\Omega_i)\cap W^s(\Omega_i))\setminus\Omega(f).
$$
In this case, there are sequences of indices $j_m,k_m\to\infty$
as $m\to\infty$ such that
$$
f^{-j_m}(p),f^{k_m}(p)\to\Omega_i,\quad m\to\infty.
$$
Since the set $\Omega_i$ is compact, we may assume that
$$
f^{-j_m}(p)\to q\in\Omega_i\;\mbox{ and}\;f^{k_m}(p)\to r\in\Omega_i.
$$
Since $\Omega_i$ contains a dense positive semi-trajectory, there exist points
$s_m\to r$ and indices $l_m>0$ such that $f^{l_m}(s_m)\to q$ as $m\to\infty$.

Clearly, if we continue the sequence 
$$
p,f(p),\dots,f^{k_m-1}(p),s_m,\dots,f^{l_m-1}(s_m),f^{-j_m}(p),\dots,f^{-1}(p)
$$
periodically with period $k_m+l_m+j_m$, we get a periodic
$d_m$-pseudotrajectory of $f$ with $d_m\to 0$ as $m\to\infty$.

Since $f\in\LPS$, there exist periodic points $p_m$ (for $m$ large enough)
such that $p_m\to p$ as $m\to\infty$, and we get the desired
contradiction with the assumption that $p\notin\Omega(f)$.
The lemma is proved.
\medskip

Lemmas 5 -- 8 show that $\LPS\subset\Omega S$.
\medskip

\section{References}

1. S. Yu. Pilyugin, {\em Shadowing in Dynamical Systems},
Lecture Notes Math., vol. 1706, Springer, Berlin, 1999.

2. K. Palmer, {\em Shadowing in Dynamical Systems. Theory and
Applications}, Kluwer, Dordrecht, 2000.

3. D. V. Anosov, {\em On a class of invariant sets of smooth
dynamical systems}, Proc. 5th Int. Conf. on Nonlin. Oscill.,
{\bf 2}, Kiev, 1970, 39-45.

4. R. Bowen, {\em Equilibrium States and the Ergodic Theory
of Anosov Diffeomorphisms}, Lecture Notes Math., vol. 470,
Springer, Berlin, 1975.

5. C. Robinson, {\em Stability theorems and hyperbolicity in
dynamical systems}, Rocky Mount. J. Math., {\bf 7}, 1977, 425-437.

6. A. Morimoto, {\em The method of pseudo-orbit tracing and
stability of dynamical systems}, Sem. Note {\bf 39},
Tokyo Univ., 1979.

7. K. Sawada, {\em Extended $f$-orbits are approximated by
orbits}, Nagoya Math. J., {\bf 79}, 1980, 33-45.

8. P. Ko\'scielniak, {\em On genericity of shadowing and
periodic shadowing property}, J. Math. Anal. Appl., {\bf 310},
2005, 188-196.

9. S. Yu. Pilyugin, {\em Variational shadowing}, Discrete Contin.
Dyn. Syst. (accepted).

10. K. Sakai, {\em Pseudo orbit tracing property and strong
transversality of diffeomorphisms of closed manifolds},
Osaka J. Math., {\bf 31}, 1994, 373-386.

11. S. Yu. Pilyugin, A. A. Rodionova, and K. Sakai,
{\em Orbital and weak shadowing properties}, Discrete Contin.
Dyn. Syst., {\bf 9}, 2003, 287-308.

12. F. Abdenur and L. J. Diaz, {\em Pseudo-orbit shadowing in the
$C^1$ topology},  Discrete Contin.
Dyn. Syst., {\bf 7}, 2003, 223-245.

13. S. Yu. Pilyugin and S. B. Tikhomirov, {\em Lipschitz shadowing
implies structural stability} (to appear).

14. S. Yu. Pilyugin, {\em Spaces of Dynamical Systems} [in Russian],
Reg. Chaotic Dynamics, Moscow-Izhevsk, 2008.

15. S. Yu. Pilyugin, K. Sakai, and O. A. Tarakanov, {\em Transversality
properties and $C^1$-open sets of diffeomorphisms with weak
shadowing}, Discrete Contin. Dyn. Syst., {\bf 9}, 2003, 287-308.

16. O. B. Plamenevskaya, {\em Weak shadowing for two-dimensional
diffeomorphisms}, Mat. Zametki, {\bf 65}, 1999, 477-480.

17. N. Aoki, {\em The set of Axiom A diffeomorphisms with no cycle}, 
Bol. Soc. Brasil. Mat. (N.S.), {\bf 23}, 1992, 21-65.

18. S. Hayashi, {\em Diffeomorphisms in $\mathcal{F}^1(M)$ satisfy Axiom A}, 
Ergod. Theory Dyn. Syst., {\bf 12}, 1992, 233-253.

19. S. Yu. Pilyugin, {\em Sets of diffeomorphisms with various limit 
shadowing properties}, J. Dynamics Differ. Equat., {\bf 19}, 2007, 747-775.

20. S. Yu. Pilyugin, {\em Introduction to Structurally Stable Systems
of Differential Equations}, Birkh\"auser-Verlag, 1994.

\end{document}